\documentclass[a4paper,12pt]{article}
\usepackage{pdfpages}
\usepackage{amsmath,amssymb, amsthm}
\usepackage{hyperref}
\usepackage{verbatim}
\usepackage{float}
\usepackage[absolute]{textpos}

\newcommand{\PP}{\mathbb{P}}
\newcommand{\EE}{\mathbb{E}}
\newcommand{\ZZ}{\mathbb{Z}}
\newcommand{\cE}{\mathcal{E}}

\theoremstyle{theorem}
\newtheorem{theorem}{Theorem}

\theoremstyle{definition}

\usepackage[top=1.25in]{geometry}

\begin{document}
\vskip1.5cm
\begin{center}
{\Large \bf The Paradox of Anti-Inductive Dice} \\
\end{center}
\vskip0.2cm
\begin{center}
Summer Eldridge$^1$, Ivo David de Oliveira$^2$, Yogev Shpilman$^3$   
\end{center}
\vskip 8pt

\begin{center}
{ \it $^1$George Mason University\\
$^2$University of Oldenburg\\
$^3$Weizmann Institute of Science} \vspace*{0.3cm}

\href{mailto:seldrid4@gmu.edu}{\tt seldrid4@gmu.edu}, \href{mailto:ivo.david@uni-oldenburg.de}{\tt ivo.david@uni-oldenburg.de}, \href{mailto:yogev.shpilman@weizmann.ac.il}{\tt yogev.shpilman@weizmann.ac.il}
\end{center}

\vglue 0.3truecm

\begin{abstract}
\noindent 
We identify a new type of paradoxical behavior in dice, where the sum of independent rolls produces a deceptive sequence of dominance relations. We call these ``anti-inductive dice". Consider a game with two players and two non-identical dice. Each rolls their die $k$ times, adding the results, and the player with the highest sum wins. For each $k$, this induces a dominance relation between dice, with $A[k]\succ B[k]$ if $A$ is more likely than $B$ to win after $k$ rolls, and vice versa.  For certain classes of dice, the limiting behavior of these relations is well-established in the literature, but the transient behavior, the subject of this paper, is less well-understood. This transient behavior, even for dice with only 4 faces, contains an immensely rich parameter space with fractal-like behavior. 
\end{abstract}

\section{Introduction }

Consider two dice, \emph{David} and \emph{Goliath}, with faces 
\begin{equation}
\begin{array}{cc}
\mathcal{D}=\left\{ 1,1,4,4,5,6\right\} ,\, & \mathcal{G}=\left\{ 0,1,2,6,6,6\right\} \end{array}
\end{equation}
respectively. This pair of dice has a curious property. Say that whoever gets the higher number wins. When you roll them once and compare the results, Goliath is more likely than David to win. The same is true if you add the result of two rolls and compare the sums: Goliath is favored. In fact, if you roll three times, five times,  six times, or more, Goliath is always more likely to win. However, if you roll exactly four times, then, and only then, David is more likely to win.

We refer to this pair of strange dice as the \emph{David vs Goliath} dice, an instance of what we call \emph{anti-inductive} dice: a new category of counterintuitive dice behavior, somewhat similar to Efron's intransitive dice \cite{R01}. In this paper, we present computational results regarding these dice, which have as of yet resisted formal statements. More specifically, we (i) formalize the notion of \emph{anti-inductive} dice, (ii) prove the opening statement concerning David vs Goliath dice, (iii) map the outcome space of 3 and 4 sided dice, (iv) identify the simplest conditions under which late inversions in the sums of rolls occur, and (v) provide compelling empirical evidence on conjectures concerning arbitrarily deceptive dice.

\section{Literature Review}

We denote an $m$-sided die with faces $d_i$ for $i = 1,..., m$ with a capital letter $D = \{d_i\}_{i=1,...,m}$. We interchangeably refer to $D$ as both the die and the probability distribution generated by sampling one of the face values with equal probability. Furthermore, we denote the random variable resulting from the sum of $k$ independent samples from  $D$  by $D[k]$. We say the die $D_1$ ``beats" $D_2$ after $k$ rolls, notated $D_1[k]\succ D_2[k]$, if 
\begin{equation}
    \PP(D_1[k]>D_2[k])>\PP(D_1[k]<D_2[k])\,.
\end{equation}
If neither beats the other, they ``tie". This induces well-defined dominance relations between dice for each value of $k$.

Dominance relations between dice are most well-known in the phenomenon of intransitive dice, that is, three dice $A,B,C$ such that $A\succ B,\; B\succ C,\; C\succ A$. For example, the \emph{magic square dice}, given by $A=\{2,6,7\}$, $B=\{1,5,9\}$ and $C=\{3,4,8\}$. Some non-transitive dice  ``reverse" the direction of the intransitivity if the sum of two rolls is compared, i.e. $A[2]\prec B[2],\; B[2]\prec C[2],\; C[2]\prec A[2]$. Compared to the ``1-sided" $Z=\{0\}$ die,  the die $D=\{-1,-1,2\}$ reverses infinitely many times, with $D[k]\succ Z[k]$ for $n\equiv 2 (\text{mod } 3)$ and $Z[k]\succ D[k]$ otherwise \cite{G97}. In fact, for any $n$, a particular set of $n$ dice can realize all possible sets of dominance relations between pairs of dice through repeated iterations \cite{S94}.

Intransitivity in dice is interesting because it violates our intuitive expectations of rational inferences. Specifically, it contradicts the following rule: ``If A beats B, and B beats C, then A beats C" or, in our notation:
\\ \\
\textbf{Rule of Transitivity:} $\ \ \ \ \ \ \ A\succ B\succ C\Rightarrow A\succ C$;\\ \\
We now discuss \textit{David} and \textit{Goliath} and our definition of ``anti-inductive dice". Formally, we say that a pair of dice $A$ and $B$ is anti-inductive if there exists an $k>0$ such that $A[i]\succ B[i],\; \forall i<k$ but $B[k]\succ A[k]$, and the dominance relation ``$\succ$" is non-periodic. Anti-inductive dice are particularly interesting for disobeying the \emph{rule of induction}, given by\\ \\
\textbf{Rule of Induction:} $\ \ \ \ \ \ \ A[i]\succ B[i]\,\forall i<k\Rightarrow A[k]\succ B[k]$;\\ \\
a foundational rule with a long philosophical history and often cited as one of the core foundations of science itself \cite{K21}. This notion generalizes any pattern of dominance relations that holds for some number of rolls.

This is meant to formally encode the ``unpredictability" of the dice. In particular, such dice depend precisely on their exact distribution even for large $n$: you can't generalize later dominance relations even with all the information collected from earlier rolls. As shown in \cite{G97}, the limiting behavior of dominance relations between dice is determined by only two summary measures (explained in Theorem \ref{the:edge}). However, the limiting behavior does not determine the die's \emph{transient behavior}, which, as we will show throughout this paper, depends richly on the particular shape of the distribution.

\section {Limiting Behavior and David vs. Goliath}

In this section, we provide the steps to prove the opening statement of this paper, namely that David defeats Goliath only when the number of rolls is four. To do this, we will first introduce the notion of Edgeworth expansions as given in \cite{BGGH18} and then point to the computer-assisted verifications needed to complete the proof of the David vs Goliath property. Full details are provided in the appendix.

Given a pair of dice $A$ and $B$, the distribution of possible differences contains within it all the information we need to understand the problem: we can restate $\PP(A[n]>B[n])$ as $\PP(A[n]-B[n]>0)$. Recall that $A+B$ is distributed as $A*B$, where $*$ denotes the convolution, i.e.
\begin{equation}
    \PP(A+B=z)=\sum_x \PP(A=x)\PP(B=z-x)\,.
\end{equation}
Therefore, if we wish to analyze a pair of dice $A$ and $B$, it is useful to define a die $\Delta=A-B$, which we refer to as the difference die, taking the form $\Delta=A-B=A*(-B)$. For David and Goliath, $\Delta = \text{Goliath} - \text{David}$, producing a difference die with $36$ faces containing all possible scores after one roll (with repetitions), which range from $-6$ to $+5$.

Convolution is commutative and associative, so instead of considering $A[n]*(-B[n])$, we consider $\Delta[n]$, and so dominance relations between iterated dice reduce to a special case of self-convolution. The central limit theorem guarantees that $\Delta [n]$ approaches a normal distribution; however, this is insufficient for the result. While the overall distribution approaches symmetry, implying the magnitude of the probability difference goes to 0, the sign of the difference remains unrestricted.

Define the tilt $T$ of a distribution $X$ as 
\begin{equation}
    T(X) \equiv \PP(X>\EE(X))-\PP(X<\EE(X))\,.
\end{equation}
In the case of nonzero third moment, the following theorem from \cite{BGGH18} completely characterizes the limiting behavior of $T(X)$ in terms of the normalized third moment, $\nu_3=\EE(X^3)/\EE(X^2)^{3/2}$ and the span, the largest integer $b$ such that all entries can be written in the form $a+bk,\; k\in \ZZ$:
\begin{theorem} \label{the:edge}
For all $n$, $$T(X[n])=\frac{L(na)}{\sqrt{2\pi n}}+\mathcal{E}$$ where $$L(c)=\frac{(-c)\text{mod } b - c\text{ mod } b }{\sigma} -\nu_3$$
and $\mathcal{E}$ is of order $o(1/\sqrt{n})$ and is defined explicitly as a function of the die (details in the appendix).
\end{theorem}

When the error term $\cE$ is small enough, the sign of the tilt is determined by the sign of $L(na)$, and if $b=1$, it depends only on the third moment. However, the convergence of this expansion is extremely slow, allowing for tens of thousands of potentially unrestrained outcomes.  We are now ready to sketch the proof for David vs. Goliath:

 We use Theorem \ref{the:edge} to prove the opening statement of David and Goliath in four steps: We first (i) calculate the span and third moment of the difference die $\Delta\equiv$ Goliath - David. We then (ii) use the explicit formula for the error $\cE$ of Theorem \ref{the:edge}. We then (iii) calculate the first $N$ such that $|L(n)|>|\cE(n)|, \forall n>N$.
Finally, we (iv) perform an exhaustive search to verify that the property holds for all values less than $N$. The combination of these four steps allows us to state our main result:

\begin{theorem}[David vs Goliath Dice]\label{the:DvsG}
    Given two dice with faces $D=\{1,1,4,4,5,6\}$ and $G=\{0,1,2,6,6,6\}$, we have that $D[k]\succ G[k]$ if and only if $k = 4$. For all other values of $k$ we have that $D[k]\prec G[k]$.
\end{theorem}

Refer to the appendix for a more detailed breakdown of this result and the proof. These particular dice were chosen for their latest double-late-inversion among integer-valued dice whose face values add up to 21, i.e., they belong to the same class as the traditional 6-sided die. All other such dice are computationally verified to invert their dominance relation after no more than 4 rolls, if they invert at all.

\section{Mapping Dice Behavior} \label{outcomeMapping}

\begin{figure}[hbt!]
    \centering
    \vspace{-1cm}
    \includegraphics[width=1\linewidth]{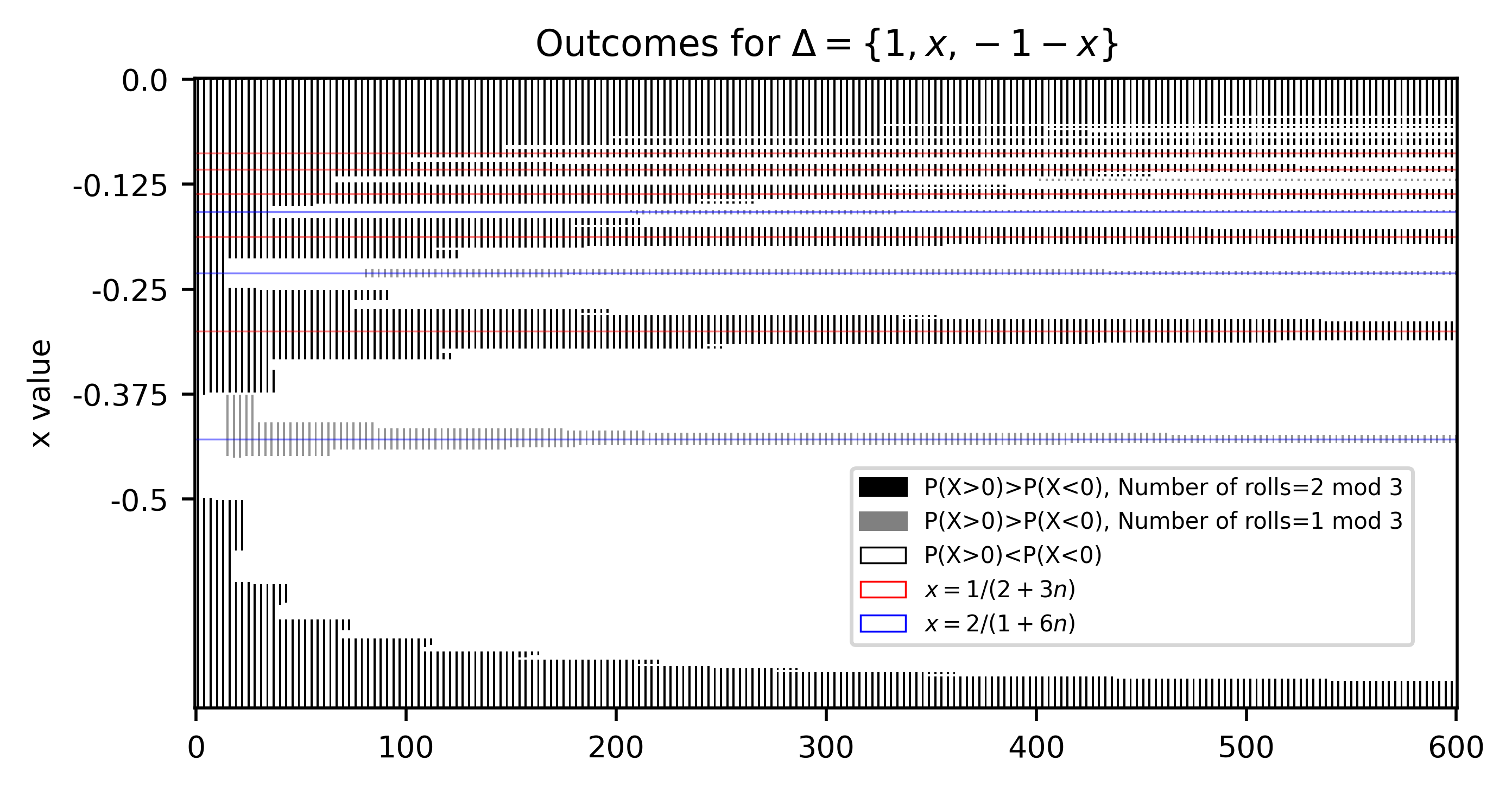}
    \vspace{-1cm}
    \caption{\small The truncated sequence of dominance relations for repeated rolls of a 3-sided die with maximum value 1, against the 0 die with one side. The centers of the ``peaks" where behavior remains periodic take a particular form, with dice that win on rolls with residue 1 and 2 mod 3 clustering around dice of the form $x=1/(2+3n)$ and dice that win on rolls with residue 1 mod 3 clustering around dice of the form $x=2/(1+6n)$.}
    \label{fig:tsided}
\end{figure}

Having proven the opening statement of this paper in Theorem \ref{the:DvsG}, we now turn our attention to the problem of (i) mapping the transient behavior of the sums of dice and more specifically we seek to (ii) understand the minimal conditions under which non-trivial anti-inductive behavior can be found. When comparing two dice, their dominance relations are unaffected by adding a constant to every face of both dice, so for the remainder of the paper, dice will be assumed to have an expected value of 0 (``balanced") unless stated otherwise. 

The space of 6-sided dice is relatively high-dimensional, too large to be dealt with directly. To render an explicit study feasible, we simplify the problem as much as possible while maintaining the original behavior. Instead of considering pairs of dice, we can consider a single die and compare it to the ``1-sided die" which always rolls 0. Two-sided balanced dice are always symmetrical, which is preserved under convolution, so they must always tie against the $\{0\}$ die, and so cannot exhibit the desired behavior.

Because 3-sided dice have only one parameter, losing one from fixing the expected value to 0 and another by normalizing the largest face to 1, the dependence of outcomes on parameters can be represented fully in a 2-dimensional plot. If the dice are normalized to assume the largest value is positive, all dice can be represented in coordinates of $\{1,x,-1-x\}$, with $0>x>-\frac{1}{2}$, as plotted in Figure \ref{fig:tsided}. With the opposite normalization, the values are exactly negative of those here, and so the outcomes have flipped sign. Note that in all observed runs, 3-sided dice never perfectly tie with the 0 die, and always lose after $3n,\, n\in \mathbb{Z}$ rolls if their largest value is positive. This corresponds to the result of Grinstead \cite{G97} that, in the limiting case, the behavior of $n$ sided die after $k=n\cdot i, i\in \mathbb{N}$ rolls depends only on the skew and the span. The skew for this half of the die-space is always negative, and in any case is continuous, as opposed to the span, which is not continuous and varies in sign even in this half-space.

\begin{figure}[H]
    \centering
    \vspace{-0.4 cm}
    \includegraphics[width=0.95\linewidth]{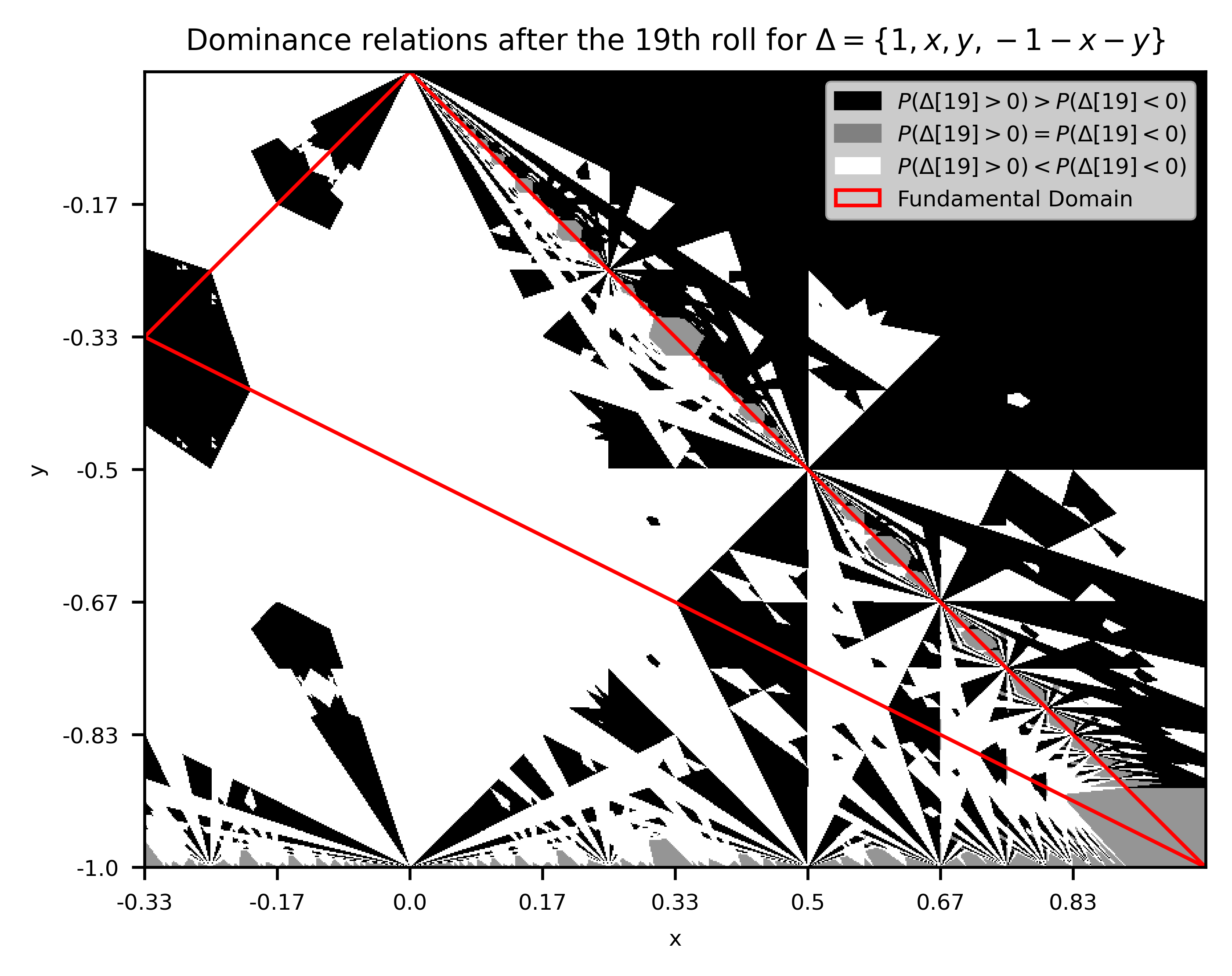}
    \vspace{-0.6 cm}
    \caption{\small Dominance relations between $\Delta[19]$ and the 0 die are plotted in the $\{1,x,y,-1-x-y\}$ parametrization. Even after less than 20 rolls, the regions where a die wins or loses against the zero dice are very elaborate. We suspect, but have not yet proved, that the limiting shape is infinitely rough.}
    \label{fig:19th}

\end{figure}

The set of all 4-sided dice with expectation value 0 is naturally 2-dimensional. We represent normalized 4-sided die in coordinates $\Delta=\{1,x,y,-1-x-y\}$, $|x|,|y|\leq 1$, but this representation is redundant. It distinguishes dice with the same faces in a different order, and it includes some dice whose largest face is negative, which have exactly inverse behavior of a particular normalized die with the largest face positive. If we remove the redundancies, assuring $x>y>-1-x-y$ and $-1-x-y>-1$, the ``fundamental domain" satisfies $-\frac{1}{3}\leq x\leq 1$ and $\frac{1-x}{2}\leq y\leq -|x|$. We are then able to depict the set of outcome sequences in a 2-dimensional map, where any observed symmetries within the domain are inherent to the space itself rather than redundant parametrization. First, we can plot the outcomes after a given number of rolls. To demonstrate the complexity exhibited, we plot the outcomes after 19 rolls in \ref{fig:19th}. We choose an odd-numbered roll, as even-numbered rolls have much simpler geometry.

To represent the way outcomes after sequential rolls influence each other, we'd like to be able to plot sequences of rolls, as in Figure \ref{fig:tsided}. Because we lack an extra dimension to represent the number of rolls, we collapse each sequence of outcomes onto a single point. Given a die $\Delta$, there is a sequence $\{B_k\}_{k=1}$, where $B_k$ encodes its dominance relation with the $\{0\}$ die, with $0,1,2$ for $\Delta[k]\prec \{0\},\; \Delta[k]\approx \{0\},\; \Delta[k]\succ \{0\}$ (``losing, tying, winning") respectively. Such a sequence has a natural representation as a number in trinary by concatenating the entries of $\{B_k\}$. This representation expresses which dice have similar outcomes, with relatively more weight placed on earlier rolls, corresponding to simpler properties. 

\begin{figure}[H]
    \centering
    \includegraphics[width=0.9\linewidth]{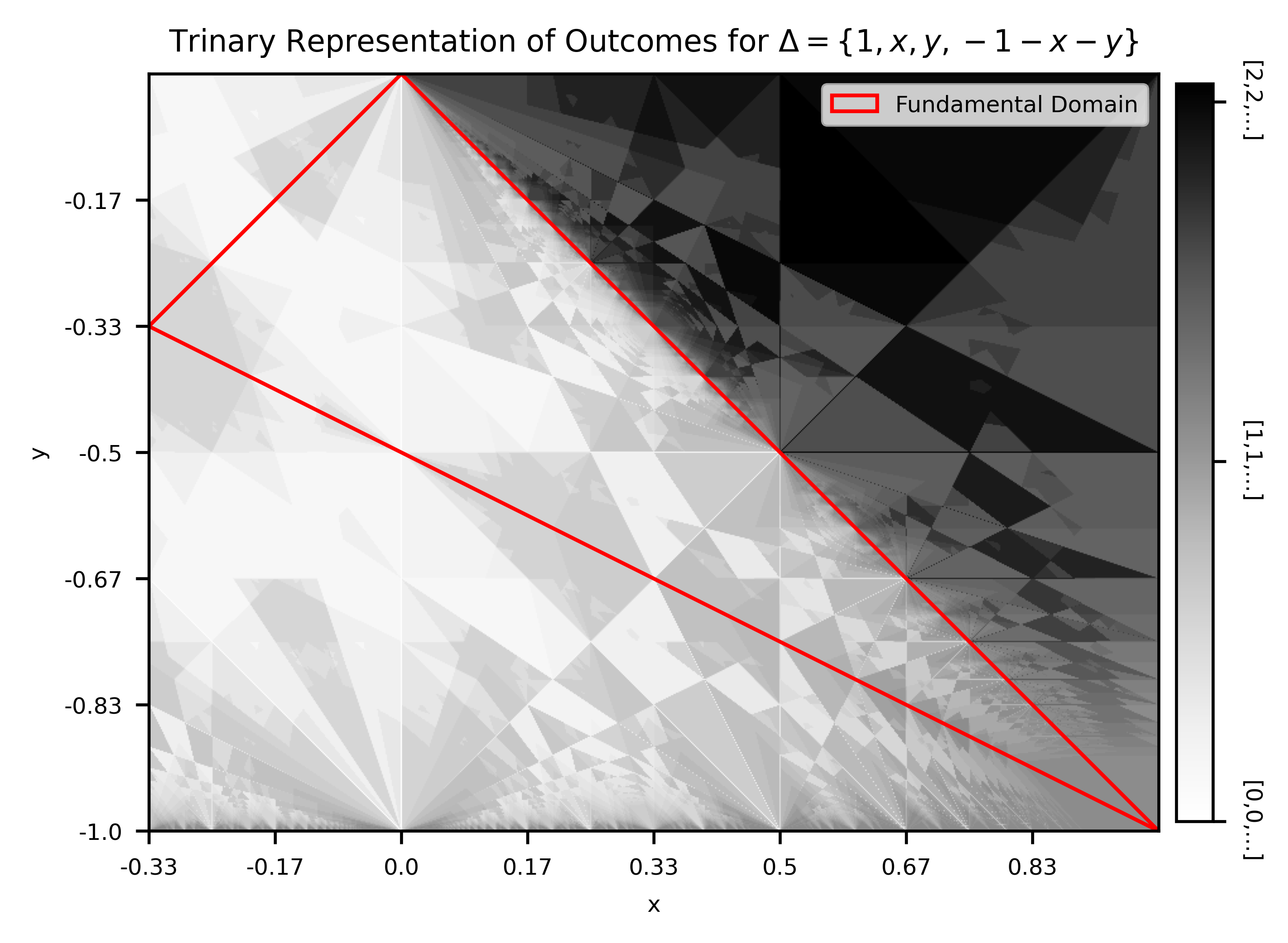}
    \caption{\small Given a die $\Delta$, there is a sequence $\{B_k\}_{k=1}$, where $B_k$ encodes its dominance relation with the $\{0\}$ die, using $0,1,2$ for $\Delta[k]\prec \{0\},\; \Delta[k]\approx \{0\},\; \Delta[k]\succ \{0\}$. The sequence is converted to a number in trinary by concatenating the values, which here are represented with shade. This expresses which areas in the parameter space have similar outcomes, with more weight placed on earlier rolls}
    \label{fig:trinar}
\end{figure}

The outcome map resulting from truncated sequences of 20 rolls is displayed in Figure \ref{fig:trinar}. Even after relatively few rolls, the complexity is immense, as thousands of interrelated equations balance out to produce a given outcome. The structure here has so far eluded formal statements. Each possible outcome is a particular linear expression of the face values, which explains the linear borders of the regions, but not their radiative structure. The centers of these ``spokes" are identified as small-denominator rational points, which will have different behavior from any close point due to exact cancellation, but why this should give them an elaborate symmetry is unclear. There is also a seeming self-similarity as you follow the fundamental domain towards its lower right corner, not perfect in this parametrization, but suggestive. Note that in the lower right corner, the uniform gray tone represents a tie for every outcome, as the die is too close to symmetrical for a distinction to appear after only 20 rolls. Only dice with rational face values were considered to allow for arbitrary precision in calculations, and so measure 0 sets along lines of perfect ratios are distinguishable where they would not be in a ``true" map. 

\section{Late Inversions}

The space of balanced and normalized 5-sided dice is 3-dimensional, so a 2-dimensional map is impossible, but we can search within the space to find particularly deceptive dice. Such a computational search has found a family of particularly deceptive dice. In particular, as explained in the section on limiting behavior of dice, if a die has a span of 1 and a positive third moment, it must have negative limiting behavior, but we have demonstrated computationally they may win for the first $k$ rolls, verified up to over $k=1000$. Purely negative limiting behavior removes the possibility of periodicity, and in fact, the result of Grinstead \cite{G97} implies no die can have limiting periodicity of minimal period greater than the die's number of sides.

\begin{figure}[hbt!]
    \centering
    \includegraphics[width=0.75\linewidth]{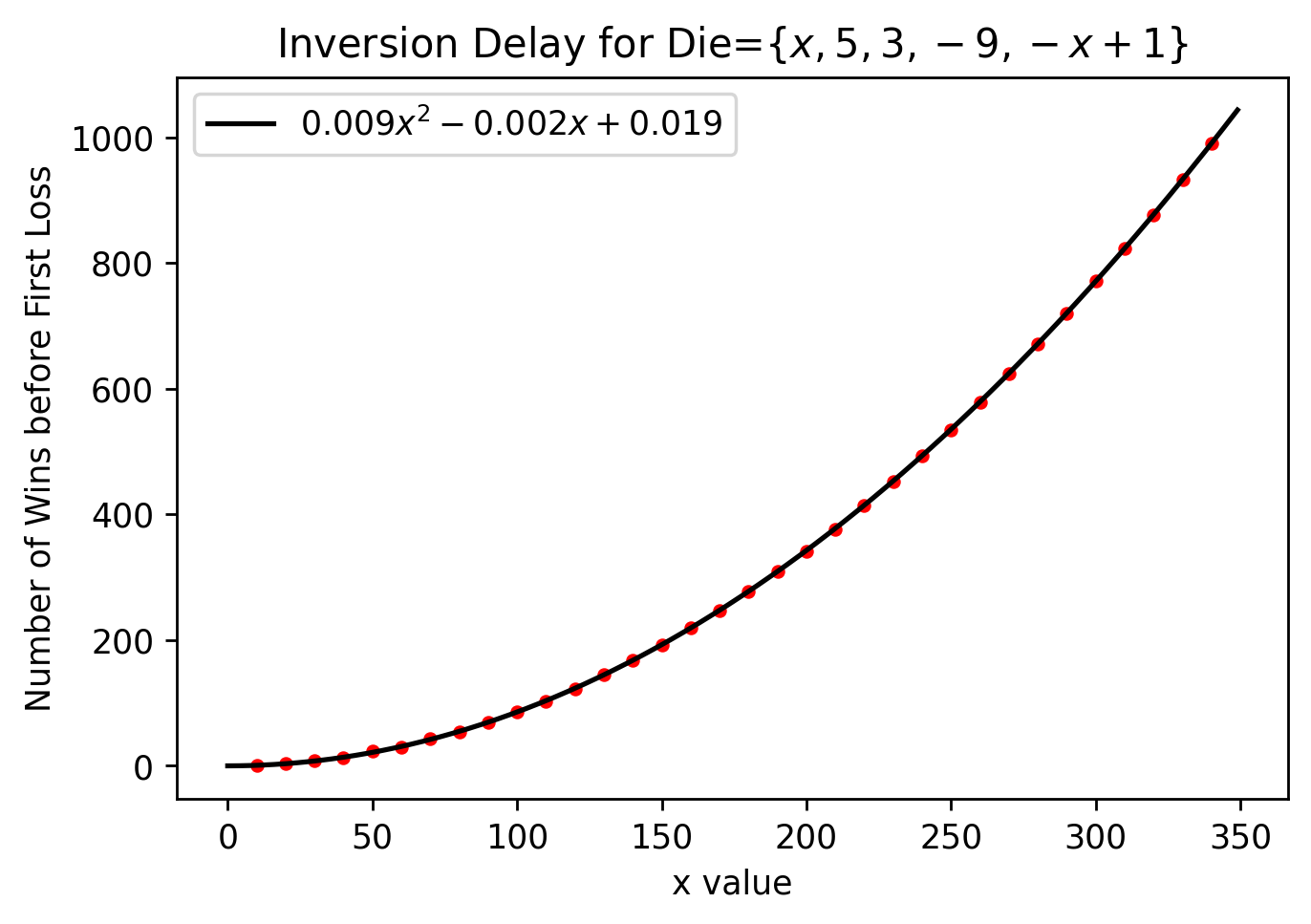}
    \caption{Approximately quadratic inversion delay, which, if it continues, would indicate a die can win for arbitrarily many rolls until a certain point, losing forever after.}
    \label{fig:invdelay}
\end{figure}

Consider the following parametrized family of dice, defined by 
\begin{equation}
    \Delta(x)=\{x,5,3,-9,-x+1\}\,.
\end{equation}
For this family, the first $k$ such that $\Delta[k]\prec \{0\}$, which we term an ``inversion", is approximately quadratic in $x$. This quadratic inversion time and the extremely late inversions exhibited by this family are plotted in Figure \ref{fig:invdelay}. We can characterize the behavior of this infinite family by breaking down the set of outcomes. Use $A$ and $B$ to refer to the number of times we roll $x$ and $-x+1$ respectively, and $k$ to refer to our total number of rolls. If $x>9k$, then for any outcome where $A>B$, the total sum will be positive, and for any outcome where $B>A$, the total sum will be negative. The two events are equally likely, so the tilt is determined entirely by the cases where $B=A$. In these cases, $x$ and $-x+1$ leave a residue of $A$, so the distribution is independent of the value of $x$. Therefore the tilt of $D(x)[k]$ is constant for all $x>9k$, and so the question of arbitrarily late inversions in this family of 5-sided dice reduces to considering the tilt of $E[k]$, which is the distribution conditional on $A=B$, independent of any particular die. $E[k]$ is notably not the result of convolution, as its mean is not linear in the number of iterations, and so its limiting behavior is not amenable to the preexisting Edgeworth expansion technique.

While we as yet lack proof that inversions can be delayed arbitrarily late, constant tilt when $k<x/9$ allows us to bound from below the time of first inversion for an infinite family of dice.
We conjecture that 5-sided dice are the smallest dice to exhibit arbitrarily late inversions in the sense above.

\section{Summary}

The transitive dominance relations between dice are full of elusive regularities. Both three and four-sided dice exhibit clear structure in their outcome diagrams, but the precise nature of this structure remains unclear. We suspect that in the case of the four-sided die, the limiting boundaries of winning and losing regions are infinitely rough, but it may be that the dynamics simplify after a sufficiently large number of rolls. Similarly, a five-sided die can exhibit very late inversions that have not been found in three or four-sided dice, and which we conjecture do not occur in them. Our formal results for the David and Goliath dice demonstrate that, for specific dice, results on extremely late behavior can be achieved through a combination of brute force calculation and results on limiting behavior, but any attempt to prove arbitrarily late inversion necessarily must handle very late behavior for an infinite class of dice, precluding purely computational approaches. We have only scratched the surface here of what is clearly a very rich area of behavior in some of the simplest distributions one can specify.

\section{Acknowledgement}
We would like to express our gratitude to Erika Roldan, Gerardo Barrera Vargas, and Tanya Khovanova for their essential conversations.

\pagebreak
\section{Appendix}

\subsection{Proof of David vs Goliath}
In what follows, we show that David never wins against Goliath after $n>58116$. The notation and calculations are all based on Theorem~18 of \cite{BGGH18}
 
Recall that $\mathcal{D}= \{1,1,4,4,5,6\}$ and $\mathcal{G}= \{0,1,2,6,6,6\}$.
Then we have $36$ differences for $\mathcal{G}-\mathcal{D}$, each one having probability $1/36$, which we can encode in the following matrix
\begin{equation}
\mathcal{G}-\mathcal{D} = 
\begin{pmatrix}
-1 & -1 & -4 & -4 & -5 & -6 \\
0 & 0 & -3 & -3 & -4 & -5 \\
1 & 1 & -2 & -2 & -3 & -4 \\
5 & 5 & 2 & 2 & 1 & 0 \\
5 & 5 & 2 & 2 & 1 & 0 \\
5 & 5 & 2 & 2 & 1 & 0 \\
\end{pmatrix}.
\end{equation}
The random variable $X$ produced by one sample of the difference die has the following property:\\ \\ $X$ has a span of $b=1$ and shift $a=0$,\\ \\ meaning that $a+bn$ for $n\in \ZZ$ contains all values of $X$, that is that $b=1$ is the GCD (greatest common divisor) of all values of $x-a$ for $x\in X$. For each $k\in \mathbb{N}$ we set $\mu_k:= \EE[X^k]$. Straightforward computations yield that
\begin{equation}   
\mu_1 = 0,\ \ \ \mu_2 = 10.166666..., \ \ \ \ \mu_3=-0.5, \ \ \ \mu_4= 213.833333....
\end{equation}
Now, we define $\sigma:= \sqrt{\mu_2}$ and $\nu_k:= \mu_k/\sigma^k$, $k\in \mathbb{N}$. Then we have
\begin{equation}  
\nu_1 = 0,\ \ \ \nu_2 = 1, \ \ \ \ \nu_3=-0.015424..., \ \ \ \nu_4= 2.068797....
\end{equation}
We set 
\begin{equation}
    \beta:= \frac{b/2 - na \mod b}{\sigma} = \frac{1/2 }{\sigma} = \frac{1}{2\sigma} = 0.156812....
\end{equation}
We then calculate the following parameters
\begin{align}
p_0:= e - 1 &\approx 1.718281...           &  p_1:= 3(\pi-3)/\pi^3 &= 0.013699...\\
q_1:= \frac{1}{5}+\frac{\nu_4}{24}&= 0.286199...         &  q_2:=\frac{p_0q_1}{2}+\frac{b^2p_1}{\mu_2}&= 0.247233...  \\
q_3:= |\beta|&= 0.156812...  &  q_4:= |\nu_3|/6&= 0.002570...      \\
q_5:= \frac{q_3^3}{6}+\frac{3q_3^2q_4}{2}+\frac{15q_3q_4^2}{2}+\frac{35q_4^3}{2}&= 0.000745...  &  r:= \frac{16b^2m}{(\pi C \sigma)^2}&= 0.001107...,             
\end{align}
where $m= 1/36$ is defined as the minimum of the probabilities for the outcome of the difference die values; and $C:=b$, defined as the $\ell^1$ norm of any vector $c$ such that $cX=b$ and $\sum_{j}c_j=0$. In our case, many values are easy to obtain: $2b$ is one example.
By Theorem~18 in \cite{BGGH18}  we have that
\begin{equation}
    \PP(X_n>0)-\PP(X_n<0) = \frac{-\nu_3}{3\sqrt{2\pi n}}+\cE,
\end{equation}
for all $n\geq \max \{q_1/4, 1/q_1, 81b^4/(q_1\pi^4\mu_2^2)\}= 3.494060...$, where the error term $\cE$ described in Theorem 1 is bounded by
\begin{equation}
    |\cE|\leq \frac{2q_2}{n}+\frac{e^{-nr/2}}{nr}+\frac{2q_5}{\sqrt{2\pi n^3}}+e^{-2\sqrt{n/q_1}}\left(\frac{1+p_0}{2\sqrt{n/q_1}} + \frac{4p_0q_1}{n} + \frac{1}{\pi\sqrt[4]{q_1n}}\left(\frac{q_3+q_4}{2\sqrt{n/q_1}}+2q_4\right) \right).
\end{equation}
Replacing the aforementioned constants, we obtain
\begin{equation}\label{eq:twoterms}
    \PP(X_n>0)-\PP(X_n<0) {=  0.002051}\frac{1}{\sqrt{n}}+\cE
\end{equation}
with
\begin{equation}
\begin{split}
   |\cE|&\leq 0.494467\frac{1}{n}+903.068802\frac{e^{-0.000553n}}{n}+0.000594\frac{1}{ n^{3/2}} \\
   &\quad
   +e^{-3.738481\sqrt{n}}\left(0.727108\frac{1}{\sqrt{n}} + 1.967088\frac{1}{n} + 0.435193\frac{1}{\sqrt[4]{n}}\left(0.042633\frac{1}{\sqrt{n}}+0.005141\right) \right).
\end{split}
\end{equation}
The first value of $n$ where the first term in the right-hand side of~\eqref{eq:twoterms} is greater than the error term is at $n = 58117$, that is,
\begin{equation}   
0.002051\frac{1}{\sqrt{n}}+\cE>0\quad \textrm{ for all }\quad n\geq 58117.
\end{equation}
A computer verification confirms that the conjecture holds for all $n$ up to $58116$.

\bibliographystyle{vancouver}
\bibliography{sources}

\begin{thebibliography}{1}

\bibitem{R01}
Rump CM.
\newblock Strategies for Rolling the Efron Dice.
\newblock Mathematics Magazine. 2001;74(3):212-6.

\bibitem{G97}
Grinstead CM.
\newblock On Medians of Lattice Distributions and a Game with Two Dice.
\newblock Combinatorics, Probability and Computing. 1997;6(3):273–294.

\bibitem{S94}
Savage RP.
\newblock The Paradox of Nontransitive Dice.
\newblock The American Mathematical Monthly. 1994;101(5):429-36.

\bibitem{K21}
Kyburg JHE.
\newblock In: Science \& Reason, Chapter 4: Induction. Oxford University Press; 2021. p. 58–73.

\bibitem{BGGH18}
Buhler JP, Gamst AC, Graham R, Hales AW.
\newblock 20.
\newblock In: Butler S, Cooper J, Hurlbert G, editors. Explicit Error Bounds for Lattice Edgeworth Expansions. Cambridge University Press; 2018. p. 321–352.

\end{thebibliography}

\end{document}